\documentclass[11pt]{amsart}

\usepackage{amsmath, amssymb, amsthm, amsfonts, enumerate, color, comment}
\usepackage{mathrsfs}
\usepackage{parskip}
\usepackage{xcolor}

\usepackage{mathtools}
\mathtoolsset{showonlyrefs}

\usepackage{palatino}
\usepackage{graphicx}

\usepackage{parskip}

\numberwithin{equation}{section}
\theoremstyle{plain}
\newtheorem{Proposition}[equation]{Proposition}
\newtheorem{Corollary}[equation]{Corollary}
\newtheorem*{Corollary*}{Corollary}
\newtheorem{Theorem}[equation]{Theorem}
\newtheorem*{Theorem*}{Theorem}

\theoremstyle{definition}

\usepackage{enumitem}
\setlist[enumerate]{leftmargin=*}
\setlist[itemize]{leftmargin=*}

\setlist[enumerate,1]{label=(\alph*),font=\upshape}

\setlist[enumerate,2]{label=(\roman*),font=\upshape}

\def\C{\mathbb{C}}

\def\D{\mathbb{D}}
\def\T{\mathbb{T}}

\usepackage{caption}
\usepackage{subcaption}

\renewcommand{\leq}{\leqslant}
\renewcommand{\geq}{\geqslant}
\renewcommand{\subset}{\subseteq}
\renewcommand{\phi}{\varphi}

 \renewcommand{\Re}[1]{\operatorname{Re} #1 }

\usepackage{xcolor}

	\author[A. Belli]{Anil Belli}
	\address{Department of Mathematics, University of Thessaloniki, 54124 Thessaloniki,
Greece}
	\email{anilbelli@math.auth.gr}
	
	\author[U. Gul]{Ugur Gul}
	\address{Hacettepe University, Department of Mathematics, 06800, Beytepe, Ankara, Turkey}
	\email{gulugur@gmail.com}

\author[W. Ross]{William T. Ross}
	\address{Department of Mathematics and Computer Science, University of Richmond, Richmond, VA 23173, USA}
	\email{wross@richmond.edu}
	
	\author[A. Siskakis]{Aristomenis G. Siskakis}
	\address{Department of Mathematics, University of Thessaloniki, 54124 Thessaloniki,
Greece}
	\email{siskakis@math.auth.gr}
	

\title{The Ces\`{a}ro operator is cyclic on $H^p$}


\begin{document}

\begin{abstract}
In this paper, we show that the Ces\`{a}ro operator on the Hardy space $H^p$, $0 < p < \infty$, is cyclic. Our techniques will involve semigroups. 
\end{abstract}

\maketitle

\section{Introduction}

Siskakis \cite{MR897683, MR1021904} showed that the classical {\em Ces\`{a}ro operator}, 
$$(Cf)(z) = \frac{1}{z} \int_{0}^{z} \frac{f(\xi)}{1 - \xi} d\xi$$ is a bounded linear operator on the {\em Hardy spaces} $H^p$, $1 \leq p < \infty$,  of the open unit disk $\D = \{z: |z| < 1\}$ (see \cite{Duren} and the discussion below for basic information about $H^p$  spaces) and identified its norm and spectral properties (see \cite{MR5037546} for an exposition of the Ces\`{a}ro operator). Subsequently,  Miao \cite{MR1104399} established the boundedness of the Ces\`{a}ro operator on $H^p$ when $0 < p < 1$. 
In this paper, we prove the following cyclicity result. In the following, let $\chi \equiv 1$ denote the constant function equal to one on the open unit disk $\D$ and $\overline{\operatorname{span}}$ denote the closed linear span in $H^p$. 

\begin{Theorem}\label{cyckucmainp}
For each $0 <  p < \infty$, the function $\chi$ satisfies 
$$\overline{\operatorname{span}}\{C^n \chi: n \geq 0\} = H^p.$$
\end{Theorem}

In the parlance of operator theory, Theorem  \ref{cyckucmainp} says that the Ces\`{a}ro operator is {\em cyclic} on $H^p$ with {\em cyclic vector }$\chi$.  When $p = 2$, this theorem is a consequence of Kriete and Trutt \cite{MR281025, MR350489} who proved that $I - C$ is unitarily equivalent to $M_z f = z f $ (multiplication by the independent variable $z$) on $P^2(\mu)$, where $P^2(\mu)$ represents the closure of the polynomials $\C[z]$ in $L^2(\mu)$ and $\mu$ is a particular finite positive measure on $\overline{\D}$. Moreover, the  isometric isomorphism $U: H^2 \to P^{2}(\mu)$ for which $U (I - C) U^{*} = M_z$ satisfies the additional property that $U \chi = \chi$. Since $\chi$ is clearly a cyclic vector for $M_z$ on $P^2(\mu)$ because the closed linear span in $L^2(\mu)$ of $\{M_{z}^n \chi: n \geq 0\}$ is $P^2(\mu)$,  one sees  that  $\chi$ is a cyclic vector for $I - C$, and hence $C$, on $H^2$.  We will see another proof of the cyclicity of $C$ when $p = 2$ in \S \ref{three} of this paper. We tried to find a proof of the cyclicity of the Ces\`{a}ro operator on $H^p$  for general $p$ but were unable to do so; hence this note. 

We point out a related paper \cite{MR4997032} in which the cyclicity of the Ces\`{a}ro operator on the (Hilbert) Hardy space of the upper-half-plane was proved using normality and Bram's cyclicity theorem. 

\section{A Reduction to the $p \geq 2$ case}\label{pl2}

For $0 <  p < \infty$, the {\em Hardy space} $H^p$ is the space of analytic functions $f$ on the open unit disk $\D$ for which 
$$\|f\|_{H^p} := \Big\{\sup_{0 < r < 1} \int_{\partial \D} |f(r \xi)|^2 dm(\xi)\Big\}^{\frac{1}{p}}$$
is finite. 
In the above, $m$ is normalized Lebesgue measure on $\partial \D$ in that $m(\partial \D) = 1$.
When $1 \leq p < \infty$, $\|\cdot\|_{H^p}$ defines a norm that makes $H^p$ a Banach space while when $0 < p < 1$, $\|f - g\|_{H^p}^{p}$ defines a metric that makes $H^p$ complete. 
By classical theory, every $f \in H^p$ has the property that the radial limit function 
$$f(\xi) = \lim_{r \to 1^{-}} f(r \xi)$$ exists and is finite for almost every $\xi \in \partial \D$. Moreover, this radial limit function, also denoted by $f$, belongs to $L^p =  L^p(\partial \D, m)$ and, in fact, $\|f\|_{L^p} = \|f\|_{H^p}$. See \cite[Ch.~2]{Duren} for the details of the above discussion. 

Since $C$ is cyclic on $H^2$ with cyclic vector $\chi$ (see the discussion in the previous section or the analysis in \S \ref{three}), the following argument will show that  $C$ is cyclic on $H^p$ for all $0 <  p < 2$ with cyclic vector $\chi$. Indeed, fix $0 <  p < 2$ and observe from H\"{o}lder's inequality, via $\|f\|_{L^p} = \|f\|_{H^p}$,  that $H^2 \subset H^p$ and 
$\|f\|_{H^p} \leq \|f\|_{H^2}$ for all $f \in H^2$. Since, as discussed in the introduction,  $\chi$ is cyclic for $C$ on $H^2$, given any $\epsilon > 0$ and $k \in \C[z]$, there is an $h \in \C[z]$ such that $$\|h(C) \chi - k\|_{H^2} < \epsilon.$$ Thus, 
$$\|h(C) \chi - k\|_{H^p} < \epsilon.$$ The known density of $\C[z]$ in $H^p$ \cite[p.~36]{Duren} shows that $\chi$ is a cyclic vector for $C$ on $H^p$. Hence, we can focus our attention on proving Theorem \ref{cyckucmainp} when $p > 2$. Along the way, we will actually reprove cyclicity in the $p = 2$ case.

\section{Semigroups and Invariant Subspaces}

Siskakis \cite{MR897683} explored the semigroup of weighted composition operators $\{S_t\}_{t \geq 0}$ on $H^p$, $1 \leq p < \infty$, defined by 
$$(S_t f)(z) = \frac{e^{-t}}{(e^{-t} - 1) z + 1}  f\Big( \frac{z e^{-t}}{(e^{-t} - 1) z + 1} \Big), \quad z \in \D, t > 0,$$ and showed \cite[Lemma 8]{MR897683} that when $p \geq 2$, the semigroup $\{S_t\}_{t \geq 0}$ is strongly continuous on $H^p$ with 
\begin{equation}\label{asd9s9d}
\|S_{t}\|_{H^p \to H^p} \leq e^{-\frac{t}{p}}, \quad t > 0.
\end{equation}
We pause for a moment to observe that for all $t \geq 0$, the function 
$$z \mapsto  \frac{z e^{-t}}{(e^{-t} - 1) z + 1}$$ is an analytic self map of $\D$ and so, by the Littlewood subordination principle \cite[p.~10]{Duren}, the composition operator 
$$f(z) \mapsto  f\Big( \frac{z e^{-t}}{(e^{-t} - 1) z + 1} \Big)$$ is bounded on $H^p$. 
The function 
$$z \mapsto  \frac{e^{-t}}{(e^{-t} - 1) z + 1}$$ is bounded and analytic on $\D$ and so each weighted composition operator $S_{t}$ a bounded on $H^p$. The norm estimate of $S_t$ in \eqref{asd9s9d} for $p \geq 2$  is the real content of the paper  \cite{MR897683} and does not hold when $1 \leq p < 2$ (although $\{S_t\}_{t \geq 0}$ still forms a strongly continuous semigroup on $H^p$ for all $p\geq 1$).

Also observed in the Siskakis paper,  based on standard facts from the theory of semigroups \cite[Ch.~10]{MR629828}, is the Laplace transform formula 
$$(\lambda I - A)^{-1} f = \int_{0}^{\infty} e^{-\lambda t} S_{t} f dt, \quad \Re \lambda > -\tfrac{1}{p}.$$
In the above, $A$ denotes the densely defined infinitesimal generator for the semigroup $\{S_t\}_{t \geq 0}$ which turns out to be 
$$A f(z) = -(1  - z) (z f(z))'$$
for $f \in  \mathcal{D}(A)$, where 
$$\mathcal{D}(A) := \{f \in H^p: (1 - z) (z f(z))' \in H^p\}$$
is the domain of $A$.
Solving a simple differential equation yields 
\begin{equation}\label{pppApppCC}
-A^{-1} = C
\end{equation}
 on $H^p$.

We now use a variation of a  trick in  Gallardo and Partington \cite{MR4757014} (also appearing in a similar analysis in \cite{CLp}). When  $p \geq 2$,  consider the related semigroup of operators 
$$\widetilde{S}_t := e^{t} S_{pt}, \quad t > 0.$$ From \eqref{asd9s9d}, this  forms a contractive semigroup in that $\|\widetilde{S}_t\| \leq 1$ for all $t > 0$. By the standard theory of semigroups,  the infinitesimal  generator $\widetilde{A}$ for $\{\widetilde{S}_t\}_{t \geq 0}$ satisfies 
$$\widetilde{A}  = p A + I.$$ Since $\{\widetilde{S}_t\}_{t \geq 0}$ is a contractive semigroup on $H^p$, $p \geq 2$, then $1$ is contained in the resolvent of $\widetilde{A}$ and thus the {\em cogenerator} $\widetilde{V}$ of $\{\widetilde{S}_t\}_{t \geq 0}$ defined by 
$$\widetilde{V} := (\widetilde{A} + I)(\widetilde{A} - I)^{-1}$$ is a bounded operator on $H^p$ and from \eqref{pppApppCC} we have 
\begin{align*}
\widetilde{V} & = (\widetilde{A} + I)(\widetilde{A} - I)^{-1}\\
& = (p A + 2 I)(p A)^{-1}\\
& = I + \tfrac{2}{p} A^{-1}\\
& = I - \tfrac{2}{p} C.
\end{align*}
Thus, $\widetilde{V}$ and $C$ on $H^p$ share the same invariant subspaces. Now use an analysis in \cite[Lemma 6.3]{CLp} (ultimately a  Banach space version of Theorem 10-9  from \cite{MR629828}) to show that the invariant subspaces for $\widetilde{V}$ and the common invariant subspaces of $\{\widetilde{S}_{t}\}_{t \geq 0}$, and hence $\{S_t\}_{t \geq 0}$,  are identical. We summarize this discussion with the following. 

\begin{Proposition}\label{sapsosdsod99}
For $p \geq 2$ and a closed subspace $\mathcal{M} \subset H^p$, the following are equivalent. 
\begin{enumerate}
\item $C \mathcal{M} \subset \mathcal{M}$; 
\item $S_t \mathcal{M} \subset \mathcal{M}$ for every $t \geq 0$. 
\end{enumerate}
\end{Proposition}

\section{Proof of the main theorem}\label{three}

By the discussion in \S \ref{pl2}, we can assume $p \geq 2$. Let $\mathcal{M}$ be the $C$-invariant subspace of $H^p$  generated by $\chi$. By Proposition \ref{sapsosdsod99}, $$S_{t} \chi  = \frac{e^{-t}}{(e^{-t} - 1) z + 1}  \in \mathcal{M}, \quad t > 0.$$ Next we observe from  \cite[p.~113]{Duren} that the dual of $H^p$ can be identified with $H^q$ ($q$ is the H\"{o}lder conjugate index to $p$) by the Cauchy dual pairing.
$$\int_{\partial \D} f \overline{g} dm, \quad f \in H^p, g \in H^q.$$
In the above,  $f$ and $g$ are understood by their almost everywhere defined radial boundary values on the unit circle  $\partial \D$ discussed earlier. 

 Thus, if $g \in H^q$ annihilates $\mathcal{M}$, then $g$ annihilates every $S_t \chi$, $t > 0$, and so 
$$\int_{\partial \D} \frac{\overline{g(\xi)}}{(e^{-t} - 1)\xi + 1} dm(\xi) = 0, \quad t > 0,$$ which, after taking complex conjugates,  is equivalent to 
$$\int_{\partial \D} \frac{1}{1 - (1 - e^{-t}) \overline{\xi}} g(\xi) dm(\xi) = 0, \quad t > 0.$$
By the Cauchy integral formula for Hardy space functions \cite[p.~40]{Duren},  the integral on the left hand side evaluates to $g(1 - e^{-t})$. This says that 
$$g(1 - e^{-t}) = 0 \; \; \mbox{for all} \; \; t > 0$$ and thus $g$ is analytic on $\D$ having a zero set with an accumulation point in $\D$. By the uniqueness theorem for analytic functions, $g \equiv 0$.  Thus, by the Hahn--Banach separation theorem, $\mathcal{M} = H^p$ and so $\chi$ is a cyclic vector for $C$ on $H^p$.

\section{A remark}

Although a  complete characterization of the cyclic vectors for the Ces\`{a}ro operator on $H^p$ is unknown, their description is a worthy problem for future study. As a humble offering of what might be possible, fix $\lambda \in \D$ and set 
$$k_{\lambda}(z) = \frac{1}{1 - \overline{\lambda} z}, \quad z \in \D.$$ 
Note that $k_{\lambda} \in H^p$ for all $\lambda \in \D$ and for all $0 < p  < \infty$.
The reader might be familiar with this function as the reproducing kernel for the Hardy space $H^2$.   As pointed out by Kriete and Trutt \cite{MR350489}, quoting an observation of Shulman, $k_{\lambda}$ is a cyclic vector for $C$ on $H^2$. 

One can check that 
$$S_{t} k_{\lambda}(z) = \frac{e^{-t}}{(e^{-t} (1 - \overline{\lambda}) - 1)z + 1}.$$
Repeating the same argument as in the previous section (assuming  $p \geq 2 $ and setting $\mathcal{M}$ to be the cyclic $C$-invariant subspace generated by $k_{\lambda}$ and observing that $S_{t} k_{\lambda} \in \mathcal{M}$ for all $t > 0$)  will show that any $g \in H^q$ that annihilates $\mathcal{M}$ must satisfy 
$$g(1 -(1 - \overline{\lambda}) e^{-t}) = 0 \; \;  \mbox{for all $t > 0$}.$$ Thus, $g$ has zeros that accumulate in $\D$ and so $g \equiv 0$. Hence, as before,  $k_{\lambda}$ is a cyclic vector for the Ces\`{a}ro operator on $H^p$ for all $p \geq 2$. In particular, the $p = 2$ case reproves the Shulman result mentioned earlier. Now use the argument from \S \ref{pl2} to conclude that $k_{\lambda}$ is a cyclic vector for the Ces\`{a}ro operator on $H^p$ for all $0 <  p < 2$. We summarize this as follows. 

\begin{Corollary}
For $0 <  p < \infty$ and $\lambda \in \D$, the function $k_{\lambda}$ is a cyclic vector for the Ces\`{a}ro operator on $H^p$. 
\end{Corollary}

\bibliographystyle{plain}

\bibliography{references}

\end{document}